\expandafter\ifx\csname mthreemacsloaded\endcsname\relax\else \fi

\magnification1100
\input amstex


 \catcode`\@=11
 \let\wlog@ld\wlog
 \def\wlog#1{\relax}

 \newif\ifIN@
 \def\m@rker{\m@@rker}
 \def\IN@{\expandafter\INN@\expandafter}
 \long\def\INN@0#1@#2@{\long\def\NI@##1#1##2##3\ENDNI@
    {\ifx\m@rker##2\IN@false\else\IN@true\fi}%
     \expandafter\NI@#2@@#1\m@rker\ENDNI@}
  \newtoks\Initialtoks@  \newtoks\Terminaltoks@
  \def\SPLIT@{\expandafter\SPLITT@\expandafter}
  \def\SPLITT@0#1@#2@{\def\TTILPS@##1#1##2@{%
     \Initialtoks@{##1}\Terminaltoks@{##2}}\expandafter\TTILPS@#2@}
  \newtoks\Trimtoks@

 \def\ForeTrim@{\expandafter\ForeTrim@@\expandafter}
 \def\ForePrim@0 #1@{\Trimtoks@{#1}}
 \def\ForeTrim@@0#1@{\IN@0\m@rker. @\m@rker.#1@%
     \ifIN@\ForePrim@0#1@%
     \else\Trimtoks@\expandafter{#1}\fi}
 
  \def\Trim@0#1@{%
      \ForeTrim@0#1@%
      \IN@0 @\the\Trimtoks@ @%
        \ifIN@
             \SPLIT@0 @\the\Trimtoks@ @\Trimtoks@\Initialtoks@
             \IN@0\the\Terminaltoks@ @ @%
                 \ifIN@
                 \else \Trimtoks@ {FigNameWithSpace}%
                 \fi
        \fi
      }

  \font\titlebold=cmbx12 scaled 1200
  \font\twelvebold=cmbx12
  \font\tenbold=cmbx10
  \font\ninebold=cmbx9
  \font\sevenbold=cmbx7
  \font\fivebold=cmbx5

  \input amssym.def \input amssym
     \font\titlemsa=msam10 at 14.4pt
     \font\titlemsb=msbm10 at 14.4pt
     \font\titleeufm=eufm10 at 14.4pt
     \font\twelvemsa=msam10 scaled 1200
     \font\twelvemsb=msbm10 scaled 1200
     \font\twelveeufm=eufm10 scaled 1200
     \font\ninemsa=msam9
     \font\ninemsb=msbm9
     \font\nineeufm=eufm9

   \ifx\cyrfam\undefined
   \else
     \immediate\write16{}%
     \message{ !!! cyr fonts already defined. !!! }
     \message{ --- edit out superfluous font defs? }
   \fi
   \newfam\cyrfam
       \font\titlecyr=wncyr10 scaled 1440 
       \font\twelvecyr=wncyr10 scaled 1200
       \font\tencyr=wncyr10
       \font\ninecyr=wncyr9
       \font\sevencyr=wncyr7
       \font\sixcyr=wncyr6

   \newfam\eusmfam
       \font\titleeusm=eusm10 scaled 1440
       \font\twelveeusm=eusm10 scaled 1200
       \font\teneusm=eusm10
       \font\nineeusm=eusm9
       \font\seveneusm=eusm7
       
       \font\fiveeusm=eusm5

\let\Cal\cal

    \font\ninemrm=cmr9 
    \font\ninei=cmmi9
    \font\ninesy=cmsy9 
    \skewchar\ninei='177
    \skewchar\ninesy='60

  \font\twelvemrm=cmr10 at 12pt 
  \font\twelvei=cmmi10 at 12pt
  \font\twelvesy=cmsy10 at 12pt

  \font\titlemrm=cmr10 at 14.4pt 
  \font\titlei=cmmi10 at 14.4pt
  \font\titlesy=cmsy10 at 14.4pt


  \def\Smallfonts{\ninepoint}

  \def\Hfont{\titlepoint\bf}
  \def\Authorfont{\twelvepoint\it}
  \def\HHfont{\twelvepoint\bf}
  \def\HHHfont{\bf}
  \def\Bibfont{\tenbf}
  \def\Coordfont{\nineit }

  \def \thfont {\bf }
  \def \pffont {\it\itSpacing }
  \def \rkfont {\bf }
  \def \dffont {\bf }
  \def \egfont {\bf }

 \def\ninepoint{%
  \def\rm{\fam0\ninerm}%
    \textfont0=\ninemrm  \scriptfont0=\sevenrm  \scriptscriptfont0=\fiverm
    \textfont1=\ninei    \scriptfont1=\seveni   \scriptscriptfont1=\fivei
  \def\mit{\fam1\ninei}%
  \def\oldstyle{\fam1\ninei}%
    \textfont2=\ninesy   \scriptfont2=\sevensy  \scriptscriptfont2=\fivesy
    \textfont3=\tenex    \scriptfont3=\tenex    \scriptscriptfont3=\tenex
  \def\it{\fam\itfam\nineit}%
    \textfont\itfam=\nineit
  \def\bf{\ifmmode\fam\bffam\else\ninebf\fi}%
    \textfont\bffam=\ninebold 
    \scriptfont\bffam=\sevenbold 
    \scriptscriptfont\bffam=\fivebold%
  \def\msa{\fam\msafam\ninemsa}%
    \textfont\msafam=\ninemsa 
    \scriptfont\msafam=\sevenmsa
    \scriptscriptfont\msafam=\fivemsa%
  \def\msb{\fam\msbfam\ninemsb}%
    \textfont\msbfam=\ninemsb%
    \scriptfont\msbfam=\sevenmsb%
    \scriptscriptfont\msbfam=\fivemsb%
  \def\eufm{\fam\eufmfam\nineeufm}%
    \textfont\eufmfam=\nineeufm
    \scriptfont\eufmfam=\seveneufm
    \scriptscriptfont\eufmfam=\fiveeufm
   \def\eusm{\fam\eusmfam\nineeusm}%
     \textfont\eusmfam=\nineeusm
     \scriptfont\eusmfam=\seveneusm
     \scriptscriptfont\eusmfam=\fiveeusm
   \def\cyr{\fam\cyrfam\ninecyr}%
     \textfont\cyrfam=\ninecyr
     \scriptfont\cyrfam=\sevencyr
     \scriptscriptfont\cyrfam=\sixcyr
  \setbox\strutbox=\hbox{\vrule
      height7pt depth3pt width0pt}%
   \baselineskip=10.8pt\rm}

 \let\eightpoint\ninepoint 

 \def\tenpoint{%
  \def\rm{\fam0\tenrm}%
    \textfont0=\tenmrm \scriptfont0=\sevenrm \scriptscriptfont0=\fiverm%
  \def\mit{\fam1\teni}%
  \def\oldstyle{\fam1\teni}%
    \textfont1=\teni   \scriptfont1=\seveni  \scriptscriptfont1=\fivei%
    \textfont2=\tensy  \scriptfont2=\sevensy \scriptscriptfont2=\fivesy%
    \textfont3=\tenex  \scriptfont3=\tenex   \scriptscriptfont3=\tenex%
  \def\it{\fam\itfam\tenit}%
    \textfont\itfam=\tenit%
  \def\bf{\ifmmode\fam\bffam\else\tenbf\fi}%
    \textfont\bffam=\tenbold
    \scriptfont\bffam=\sevenbold%
    \scriptscriptfont\bffam=\fivebold%
  \def\msa{\fam\msafam\tenmsa}%
    \textfont\msafam=\tenmsa%
    \scriptfont\msafam=\sevenmsa%
    \scriptscriptfont\msafam=\fivemsa%
  \def\msb{\fam\msbfam\tenmsb}%
    \textfont\msbfam=\tenmsb%
    \scriptfont\msbfam=\sevenmsb%
    \scriptscriptfont\msbfam=\fivemsb%
  \def\eufm{\fam\eufmfam\teneufm}%
   \textfont\eufmfam=\teneufm
   \scriptfont\eufmfam=\seveneufm
   \scriptscriptfont\eufmfam=\fiveeufm
   \def\eusm{\fam\eusmfam\teneusm}%
    \textfont\eusmfam=\teneusm
    \scriptfont\eusmfam=\seveneusm
    \scriptscriptfont\eusmfam=\fiveeusm
   \def\cyr{\fam\cyrfam\tencyr}%
    \textfont\cyrfam=\tencyr
    \scriptfont\cyrfam=\sevencyr
    \scriptscriptfont\cyrfam=\sixcyr
  \setbox\strutbox=\hbox{\vrule %
      height8.5pt depth3.5ptwidth0pt}%
  \baselineskip=\StdBaselineskip\rm}

 \def\twelvepoint{%
  \def\rm{\fam0\twelverm}%
    \textfont0=\twelvemrm \scriptfont0=\tenmrm \scriptscriptfont0=\sevenrm
    \textfont1=\twelvei   \scriptfont1=\teni   \scriptscriptfont1=\seveni
  \def\mit{\fam1\twelvei}%
  \def\oldstyle{\fam1\twelvei}%
    \textfont2=\twelvesy  \scriptfont2=\tensy  \scriptscriptfont2=\sevensy
    \textfont3=\tenex  \scriptfont3=\tenex  \scriptscriptfont3=\tenex
  \def\it{\fam\itfam\twelveit}%
    \textfont\itfam=\twelveit
  \def\bf{\ifmmode\fam\bffam\else\twelvebf\fi}%
    \textfont\bffam=\twelvebold
    \scriptfont\bffam=\tenbold%
    \scriptscriptfont\bffam=\sevenbold%
  \def\msa{\fam\msafam\twelvemsa}%
    \textfont\msafam=\twelvemsa%
    \scriptfont\msafam=\tenmsa%
    \scriptscriptfont\msafam=\sevenmsa%
  \def\msb{\fam\msbfam\twelvemsb}%
    \textfont\msbfam=\twelvemsb%
    \scriptfont\msbfam=\tenmsb%
    \scriptscriptfont\msbfam=\sevenmsb%
  \def\eufm{\fam\eufmfam\twelveeufm}%
   \textfont\eufmfam=\twelveeufm
   \scriptfont\eufmfam=\teneufm
   \scriptscriptfont\eufmfam=\seveneufm
   \def\eusm{\fam\eusmfam\twelveeusm}%
    \textfont\eusmfam=\twelveeusm
    \scriptfont\eusmfam=\teneusm
    \scriptscriptfont\eusmfam=\seveneusm
   \def\cyr{\fam\cyrfam\tencyr}%
    \textfont\cyrfam=\twelvecyr
    \scriptfont\cyrfam=\tencyr
    \scriptscriptfont\cyrfam=\sevencyr
  \setbox\strutbox=\hbox{\vrule
      height10.2pt depth4.55pt width0pt}%
  \baselineskip=14pt\rm}

 \def\titlepoint{%
    \textfont0=\titlemrm \scriptfont0=\twelvemrm \scriptscriptfont0=\tenmrm
    \textfont1=\titlei   \scriptfont1=\twelvei   \scriptscriptfont1=\teni
  \def\mit{\fam1\titlei}%
  \def\oldstyle{\fam1\titlei}%
    \textfont2=\titlesy  \scriptfont2=\twelvesy  \scriptscriptfont2=\tensy
    \textfont3=\tenex
    \scriptfont3=\tenex
    \scriptscriptfont3=\tenex
  \def\it{\fam\itfam\titleit}%
    \textfont\itfam=\titleit
  \def\bf{\ifmmode\fam\bffam\else\titlebf\fi}%
    \textfont\bffam=\titlebold
    \scriptfont\bffam=\twelvebold%
    \scriptscriptfont\bffam=\tenbold%
  \def\msa{\fam\msafam\titlemsa}%
    \textfont\msafam=\titlemsa%
    \scriptfont\msafam=\twelvemsa%
    \scriptscriptfont\msafam=\tenmsa%
  \def\msb{\fam\msbfam\titlemsb}%
    \textfont\msbfam=\titlemsb%
    \scriptfont\msbfam=\twelvemsb%
    \scriptscriptfont\msbfam=\tenmsb%
  \def\eufm{\fam\eufmfam\titleeufm}%
    \textfont\eufmfam=\titleeufm
    \scriptfont\eufmfam=\twelveeufm
    \scriptscriptfont\eufmfam=\teneufm
   \def\eusm{\fam\eusmfam\titleeusm}%
     \textfont\eusmfam=\titleeusm
     \scriptfont\eusmfam=\twelveeusm
     \scriptscriptfont\eusmfam=\teneusm
   \def\cyr{\fam\cyrfam\tencyr}%
    \textfont\cyrfam=\titlecyr
    \scriptfont\cyrfam=\twelvecyr
    \scriptscriptfont\cyrfam=\tencyr
  \setbox\strutbox=\hbox{\vrule
      height12.3pt depth5.54pt width0pt}%
  \baselineskip=16pt\rm}

\newbox\AuthorBox\newbox\TitleBox
\newbox\TFLinebox
\newbox\FLinebox
\newbox\HLinebox
\def\SetTFLinebox#1{\setbox\TFLinebox=\hbox{#1}}
\def\SetFLinebox#1{\setbox\FLinebox=\hbox{#1}}
\def\SetHLinebox#1{\setbox\HLinebox=\hbox{#1}}

 \def\SetAuthorHead#1{%
     \setbox\AuthorBox=\hbox{\ninepoint \it 
           \ignorespaces\frenchspacing#1\unskip}}
 \def\SetTitleHead#1{%
     \setbox\TitleBox=\hbox{\ninepoint \it
           \ignorespaces\frenchspacing#1\unskip}}

  \def\itSpacing{\relax}
  \def\itSpacingOff{\relax}


 \def\Hrule{\hrule width0pt height0pt}

  \newskip\ProcSkip \ProcSkip 8pt plus2pt minus2pt

 \newskip\LastSkip
 \def\SaveLastSkip{\LastSkip\lastskip}
 \def\RestoreLastSkip{\vskip-\LastSkip\vskip\LastSkip}

 \def\NoindentAfter{\everypar={\setbox0=\lastbox\everypar={}}}

 \long\def\H#1\par#2\par{\notenumber=0 \titlepagetrue%
    {
    \baselineskip=20pt
    \parindent=0pt\parskip=0pt\frenchspacing
    \leftskip=0pt plus .2\hsize minus .3\hsize
    \rightskip=0pt plus .2\hsize minus .3\hsize
 \def\\{\unskip\break}%
    \pretolerance=10000 \Hfont #1\unskip\break
     \vskip7pt\Hrule
\hfill \Authorfont #2\hfill\hfill\unskip}
    \vskip48pt plus 4pt minus 4pt
    \par\NoindentAfter\rm}

 \long\def\Hi#1\par#2\par{\notenumber=0 \titlepagetrue%
    {  \baselineskip=0pt  \parindent=0pt\parskip=0pt\frenchspacing
    \leftskip=0pt plus .2\hsize minus .3\hsize
    \rightskip=0pt plus .2\hsize minus .3\hsize
}
    \rm}


 \newdimen\PageRemainder
  \def\SetPageRemainder{
     \PageRemainder=\pagegoal
     \ifdim\PageRemainder=\maxdimen\PageRemainder=\vsize
     \else\advance\PageRemainder by -1\pagetotal\fi}

  \def\Rpt@{}\def\Rpt@@{}

  \long\def\HH#1\par{\par
  \SaveLastSkip\removelastskip\goodbreak
  \ifdim\LastSkip<30pt 
     \LastSkip 30pt
plus 3pt minus 2pt\fi
  \SetPageRemainder\advance\PageRemainder-\LastSkip
  \ifdim\PageRemainder<150pt
       \edef\Rpt@{remain = \the\PageRemainder\noexpand\\
                pagetotal=\the\pagetotal\noexpand\\
                           pagegoal=\the\pagegoal}%
          \fi
   \ifdim\PageRemainder<65pt 
       \ifdim\PageRemainder > 0pt
          \edef\Rpt@@{\noexpand\\
                      Had HH PageRemainder$<$\relax 65pt\noexpand\\
                      Hence forced break!}%
     \vskip 0pt plus .2\PageRemainder\eject 
    \fi\fi
    \vskip\LastSkip\Hrule 
    \pretolerance=10000\rightskip=0pt plus 3em
    \hangafter1 \hangindent=2.2em%
    \noindent
    \HHfont \unskip \Ednote{\Rpt@\Rpt@@}%
            \def\Rpt@{}\def\Rpt@@{}%
            \ignorespaces
            #1\par\rightskip=0pt\pretolerance=\StdPretolerance%
    \NoindentAfter
\tenpoint\rm%
     \medskip \vskip\ProcSkip}

  \long\def\HHH#1\par{\par%
  \SaveLastSkip\removelastskip\goodbreak
  \ifdim\LastSkip<\ProcSkip%
     \LastSkip\ProcSkip\fi
  \SetPageRemainder\advance\PageRemainder-\LastSkip
  \ifdim\PageRemainder<150pt
       \edef\Rpt@{remain = \the\PageRemainder\noexpand\\
                pagetotal=\the\pagetotal\noexpand\\
                           pagegoal=\the\pagegoal}%
       \fi
   \ifdim\PageRemainder<48pt  
        \ifdim\PageRemainder > 0pt
             \edef\Rpt@@{\noexpand\\
                      Had HHH PageRemainder$<$\relax48pt\noexpand\\
                      Hence forced break!}%
       \vskip 0pt plus .2\PageRemainder\eject 
      \fi\fi
   \vskip\LastSkip\par\noindent
   \HHHfont \unskip\Ednote{\Rpt@\Rpt@@}%
  \def\Rpt@{}\def\Rpt@@{}%
  \ignorespaces
   #1\unskip.\quad\rm\ignorespaces
   \ignorepars}

  \long\def\ignorepars#1\par{\def\Test{#1}%
     \ifx\Test\Empty\def\This{\ignorepars}%
        \else\def\This{\Test\par}\fi
           \This}
  \def\Empty{}

 \def\Abstract#1\par{\bgroup\Smallfonts\narrower\HHH #1\par}
 \def\endAbstract{\par\egroup}


 \def\ProcBreak{\par%
    \ifdim\lastskip<8pt%
    \removelastskip%
    \penalty-200\vskip\ProcSkip\fi}

 \def\th#1\par{\ProcBreak \noindent
   {\thfont\ignorespaces
    #1\unskip.}\it\itSpacing\kern.4em\ignorepars}

 \def\endth{\ProcBreak\rm\itSpacingOff }


 \def\pf#1\par{\ProcBreak %
    \noindent\pffont#1\unskip.\rm\itSpacingOff{\kern .7em}\ignorepars}

 \def\endpf{\medskip \ProcBreak } 

  \def\qedbox{\hbox{\vbox{
    \hrule width0.2cm height0.2pt
    \hbox to 0.2cm{\vrule height 0.2cm width 0.2pt
             \hfil\vrule height0.2cm width 0.2pt}
    \hrule width0.2cm height 0.2pt}\kern1pt}}

  \def\qed{\ifmmode\qedbox
    \else\unskip\ \hglue0mm\hfill\qedbox\ProcBreak\fi}

  \def \rk #1\par{\ProcBreak
     \noindent{\rkfont\ignorespaces #1\unskip.}%
     \rm\kern.6em\ignorepars}

  \def \endrk {\medskip\ProcBreak }

  \def \df #1\par{\ProcBreak
     \noindent{\dffont\unskip\ignorespaces #1\unskip.}%
     \rm\kern.6em\ignorepars}

  \def \enddf {\medskip\ProcBreak }

  \def \eg #1\par{\ProcBreak
     \noindent\egfont\unskip\ignorespaces #1\unskip.
     \rm\kern.6em\ignorepars}

  \def \endeg {\medskip\ProcBreak }

  \newdimen\Overhang

   \def\MaxTag@#1#2#3#4#5{\setbox0=\hbox{#4\ignorespaces#2\unskip}%
     \dimen0=\wd0\advance\dimen0 by#3
     \ifdim\dimen0<#5\relax\dimen0=#5\fi
     \expandafter\edef\csname #1Hang\endcsname{\the\dimen0}}

 \def\MaxItemTag#1{\MaxTag@{Item}{#1}{.4em}{\ItemStyle}{\parindent}}%
 \def\MaxItemItemTag#1{%
        \MaxTag@{ItemItem}{#1}{.4em}{\ItemItemStyle}{\parindent}}
 \def\MaxNrTag#1{\MaxTag@{Nr}{#1}{.5em}{\NrStyle}{\parindent}}
 \def\MaxReferenceTag#1{%
        \MaxTag@{Reference}{[#1]}{.6em}{\ninerm}{\parindent}}
 \def\MaxFootTag#1{\MaxTag@{Foot}{#1}{.4em}{\ninerm}{\z@}}

  \def\SetOverhang@{\Overhang=.8\dimen0%
     \advance\Overhang by \wd0\relax
     \ifdim\Overhang>\hangindent\relax
       \advance\Overhang by .25\dimen0%
       \Ednote{Tag is pushing text.}\osumess{Tag is pushing text.}%
     \else\Overhang=\hangindent
     \fi}

   \def\Item#1{\par\noindent
      \hangafter1\hangindent=\ItemHang
      \setbox0=\hbox{\ItemStyle\ignorespaces#1\unskip}%
      \dimen0=.4em\SetOverhang@
      \rlap{\box0}\kern\Overhang\ignorespaces}

   \def\ItemItem#1{\par\noindent
      \hangafter1\hangindent=\ItemItemHang
      \setbox0=\hbox{\ItemItemStyle\ignorespaces#1\unskip}%
      \dimen0=.4em\SetOverhang@
      \advance\hangindent by \ItemHang
      \kern\ItemHang\rlap{\box0}%
      \kern\Overhang\ignorespaces}

  \def\Nr#1{\par\noindent\hangindent=\NrHang 
    \setbox0=\hbox{\NrStyle\ignorespaces#1\unskip}%
    \dimen0=.5em\SetOverhang@
    \rlap{\box0}\kern\Overhang
    \hangindent=\z@\ignorespaces}

   \newskip\Rosterskip\Rosterskip 1pt plus1pt 
   \def\Roster{\par\ifdim\lastskip<\Rosterskip\removelastskip\vskip\Rosterskip\fi
    \bgroup}
   \def\endRoster{\par\global\edef\LastSkip@{\the\lastskip}\removelastskip
       \egroup\penalty-50\LastSkip\LastSkip@\relax
       \ifdim\LastSkip<\Rosterskip\LastSkip\Rosterskip\fi
       \vskip\LastSkip}




 \def\cite#1{
    \def\nextiii@##1,##2\end@{{\frenchspacing\rm 
      \lBr\ignorespaces##1\unskip{\rm,~\ignorespaces##2}\rBr}}%
    \IN@0,@#1@%
    \ifIN@\def\next{\nextiii@#1\end@}\else
    \def\next{{\rm\lBr#1\rBr}}\fi\next}


   \def \Bib#1\par{%
       \par\removelastskip\SetPageRemainder
       \ifdim\PageRemainder < 97pt
        \ifdim\PageRemainder > 0pt
        \vfill\eject
       \fi\fi
    \ProcBreak \par\begingroup\parskip=0 pt%
    \goodbreak \vskip 15 pt plus 10 pt
    \noindent\null\hfill\Bibfont
      \ignorespaces #1\unskip\hfill\null\par 
    \frenchspacing \Smallfonts\rm
    \parskip=2.5 pt plus 1 pt minus.5pt%
    \nobreak\vskip 12pt plus 2pt minus2pt\nobreak
    \leftskip=0 pt \baselineskip=10.5pt}

 \def\ReferenceTagSlide{0em}
  \def\ReferenceTagGap{.5em}

  \def \rf#1{\par\noindent
     \hangafter1\hangindent=\ReferenceHang      
     \setbox0=\hbox{\ninerm[\ignorespaces#1\unskip]}%
     \dimen0=\ReferenceTagGap\SetOverhang@
     \rlap{\kern\ReferenceTagSlide\box0}%
     \kern\Overhang\ignorespaces}

  \def\ref#1\par#2\par#3\par#4\par{%
     \rf{#1}#2\unskip,\ #3\unskip,\
     #4\unskip.}

  \def\endBib{\par\endgroup\vskip 12pt minus 6pt }


  \long\def\Coordinates#1\endCoordinates{
 {\par\vskip4pt\def\\{\unskip, }\Coordfont\baselineskip10.5pt\noindent#1}}

 \def\pagecontents{
  \gdef\Pagetot@l{\pagetotal}
  \ifvoid\TRMargIns\else
    \rlap{\kern\hsize\kern10pt\vbox to 0pt{%
         \box\TRMargIns\vss}}\fi
  \ifvoid\topins\else\unvbox\topins\fi
   \dimen@=\dp\@cclv \unvbox\@cclv 
   \ifvoid\footins\else 
     \vskip\skip\footins
     \footnoterule
     \unvbox\footins\fi
   \ifr@ggedbottom \kern-\dimen@ \vfil \fi}


 \newcount\Ht 

 \def \Acc{\expandafter } 

 \def\swthat{\raise -1.1 ex\hbox{\sam$\widehat{}$}}
 \def\swttilde{\raise -1.2 ex\hbox{\sam$\widetilde{}$}}
 \def \overdot{{\raise .2 ex \hbox to 0pt {\hss\bf\smash{.}\hss}}}
 \def \overcircle{{\raise .1 ex \hbox to 0pt
    {\sam$\eightpoint\scriptstyle\hss\circ\hss$}}}

 \def \Mathaccent#1#2{{\sam 
  \setbox4=\hbox{$\vphantom{#2}$}
  \Ht=\ht4 
  \setbox5=\hbox{${#1}$}
  \setbox6=\hbox{${#2}$}
  \setbox7=\hbox to .5\wd6{}
  \copy7\kern .1\Ht \raise\Ht sp\hbox{\copy5}\kern-.1\Ht
  \copy7\llap{\box6}
  }}

  \def\SwtCheck #1{
        \ifmmode \check{#1}%
                \else \v {#1}%
                \fi}

 \def\barpartial {%
   \kern .17 em
    \overline {\kern -.17 em\partial\kern-.03 em}%
    \kern .03 em}

 
  \def\Overline#1{\setbox1=\hbox{\sam ${#1}$}%
      \ifdim \wd1 > 6pt
    \kern .11 em
    \overline {\kern -.11 em#1\kern-.14 em}
    \kern .14 em
  \else
    \kern .03 em
    \overline {\kern -.03 em#1\kern-.04 em}
    \kern .04 em
  \fi}

 \def\SOverline#1{\setbox1=\hbox{\sam ${#1}$}%
      \ifdim \wd1 > 7pt
    \kern .22 em
    \overline {\kern -.22 em#1\kern-.09 em}%
    \kern .09 em
  \else
    \kern .10 em
    \overline {\kern -.10 em#1\kern-.04 em}%
    \kern .04 em
  \fi}


 \def\Underline#1{\setbox1=\hbox{\sam ${#1}$}%
      \ifdim \wd1 > 6pt
    \kern .11 em
    \underline {\kern -.11 em#1\kern-.14 em}
    \kern .14 em
  \else
    \kern .03 em
    \underline {\kern -.03 em#1\kern-.04 em}
    \kern .04 em
  \fi}

 \def\SUnderline#1{\setbox1=\hbox{\sam ${#1}$}%
      \ifdim \wd1 > 7pt
    \kern .04 em
    \underline {\kern -.04 em#1\kern-.2 em}%
    \kern .2 em
  \else
    \kern .0 em
    \underline {\kern -.0 em#1\kern-.15 em}%
    \kern .15 em
  \fi}


 \def \Blackbox
   {\leavevmode\hskip .3pt \vbox
   {\hrule height 5pt\hbox{\hskip 4.5pt}}\hskip .5pt}

 \def \XX{\Blackbox\kern.5pt\Blackbox} 

  \def\.{.\kern1pt}

    \def\Hyphen{\edef\this{\the\hyphenchar\font}%
          \hyphenchar\font=-1\char\this\hyphenchar\font=\this}

 \ifx\undefined\text
  \def\text#1{\hbox{\rm #1}}\fi 



   \everymath{}  

  \def\PassMath@@{\aftergroup\AfterMath@} 

 \let\PassMath@\PassMath@@

 \def\AfterMath@{\futurelet\next\AfterMathMole@}

 \def\AfterMathMole@{
      \ifcat\next\space
          \def\this{}
      \else
      \ifcat\next\egroup %
        \def\this{\osumess{Handset mathsurround?? ---(see dollar brace)}}%
      \else
      \def\this{\AAfterMath@}
      \fi\fi
      \this}

 \def\hyphen@{-}
 \def\paren@{)}
 \def\apostr@{'}

 \def\MSC#1{\kern-.8\mathsurround#1\kern.8\mathsurround}

 \def\AAfterMath@#1{\def\Next{#1}
    \IN@0\Next @,.;:!?\relax @%
    \ifIN@\def\this{\MSC{\Next}}%
    \else
    \ifx\Next\hyphen@\def\this{\futurelet\next\AfterHyphen@}%
    \else
    \ifx\Next\paren@\def\this{#1}%
    \else 
    \ifx\Next\apostr@\def\this{#1}%
    \else \def\this{\osumess{Handset mathsurround??}%
                 #1}\fi\fi\fi\fi
    \this}

 \def\AfterHyphen@#1{\def\Next{#1}%
   \ifx\Next\hyphen@\def\this{--}\else
   \ifcat\next\space%
   \def\this{\kern-\mathsurround\kern.05em- \Next}\else
   \def\this{\kern-\mathsurround\kern.05em\Hyphen\Next}\fi\fi\this}

 \def\sam{\mathsurround=\z@\let\PassMath@\relax}  %
 \def\mas{\mathsurround=\StdMathsurround\let\PassMath@\PassMath@@}
 
 \def\Mas{\mathsurround=\StdMathsurround
                \everymath{\PassMath@}\let\PassMath@\PassMath@@}

 \def\m@th{\mathsurround=\z@\everymath{}}

 \def\m@@th{\mathsurround=\z@\everymath={}\let\m@th\relax}

\def\underbar#1{$\setbox\z@\hbox{#1}\dp\z@\z@
      \m@th \underline{\box\z@}$\relax}

\def\mathhexbox#1#2#3{\leavevmode
  \hbox{\m@@th$\m@th \mathchar"#1#2#3$}}

\def\dots{\relax\ifmmode\ldots\else$\m@th\ldots\,$\relax\fi}

\def\dotfill{\cleaders\hbox{\m@@th$\m@th \mkern1.5mu.\mkern1.5mu$}\hfill}
\def\rightarrowfill{$\m@th\mathord-\mkern-6mu%
  \cleaders\hbox{\m@@th$\mkern-2mu\mathord-\mkern-2mu$}\hfill
  \mkern-6mu\mathord\rightarrow$\relax}
\def\leftarrowfill{$\m@th\mathord\leftarrow\mkern-6mu%
  \cleaders\hbox{\m@@th$\mkern-2mu\mathord-\mkern-2mu$}\hfill
  \mkern-6mu\mathord-$\relax}

\def\downbracefill{$\m@th\braceld\leaders\vrule\hfill\braceru
  \bracelu\leaders\vrule\hfill\bracerd$\relax}
\def\upbracefill{$\m@th\bracelu\leaders\vrule\hfill\bracerd
  \braceld\leaders\vrule\hfill\braceru$\relax}

\def\angle{{\vbox{\m@@th\ialign{$\m@th\scriptstyle##$\crcr
      \not\mathrel{\mkern14mu}\crcr
      \noalign{\nointerlineskip}
      \mkern2.5mu\leaders\hrule height.34pt\hfill\mkern2.5mu\crcr}}}}

\def\big#1{{\m@@th\hbox{$\left#1\vbox to8.5\p@{}\right.\n@space$}}}
\def\Big#1{{\m@@th\hbox{$\left#1\vbox to11.5\p@{}\right.\n@space$}}}
\def\bigg#1{{\m@@th\hbox{$\left#1\vbox to14.5\p@{}\right.\n@space$}}}
\def\Bigg#1{{\m@@th\hbox{$\left#1\vbox to17.5\p@{}\right.\n@space$}}}
\def\n@space{\nulldelimiterspace\z@ \m@th}

\def\root#1\of{\setbox\rootbox\hbox{\m@@th$\m@th\scriptscriptstyle{#1}$}
  \mathpalette\r@@t}
\def\r@@t#1#2{\setbox\z@\hbox{\m@@th$\m@th#1\sqrt{#2}$\relax}
  \dimen@\ht\z@ \advance\dimen@-\dp\z@
  \mkern5mu\raise.6\dimen@\copy\rootbox \mkern-10mu \box\z@}

\def\mathph@nt#1#2{\setbox\z@\hbox{\m@@th$\m@th#1{#2}$}\finph@nt}

\def\mathsm@sh#1#2{\setbox\z@\hbox{\m@@th$\m@th#1{#2}$}\finsm@sh}

\def\@vereq#1#2{\lower.5\p@\vbox{\m@@th\baselineskip\z@skip\lineskip-.5\p@
    \ialign{$\m@th#1\hfil##\hfil$\crcr#2\crcr=\crcr}}}

\def\mathpalette#1#2{\sam\mathchoice{#1\displaystyle{#2}}%
  {#1\textstyle{#2}}{#1\scriptstyle{#2}}{#1\scriptscriptstyle{#2}}\mas}

\def\widehat#1{\setbox\z@\hbox{\sam$#1$}%
 \ifdim\wd\z@>\tw@ em\mathaccent"0\msbfam@5B{#1}%
 \else\mathaccent"0362{#1}\fi}
\def\widetilde#1{\setbox\z@\hbox{\sam$#1$}%
 \ifdim\wd\z@>\tw@ em\mathaccent"0\msbfam@5D{#1}%
 \else\mathaccent"0365{#1}\fi}

 \def\dots{\relax{}
  \ifmmode\def\thedots{\mdots@}\else\def\thedots{\tdots@}\fi %
  \thedots}

 \let\@ldeqno\eqno\let\@ldleqno\leqno
 \def\eqno{\everymath{}\@ldeqno} \def\leqno{\everymath{}\@ldleqno}

  \let\@ldeqalignno\eqalignno
  \def\eqalignno#1{\sam\@ldeqalignno{#1}\mas}
  \let\@ldeqalign\eqalign
  \def\eqalign#1{\sam\@ldeqalign{#1}\mas}

 \def\overrightarrow#1{\vbox{\m@th\ialign{##\crcr
      \rightarrowfill\crcr\noalign{\kern-\p@\nointerlineskip}
      $\hfil\displaystyle{#1}\hfil$\crcr}}}
 \def\overleftarrow#1{\vbox{\m@th\ialign{##\crcr
      \leftarrowfill\crcr\noalign{\kern-\p@\nointerlineskip}
      $\hfil\displaystyle{#1}\hfil$\crcr}}}
 \def\overbrace#1{\mathop{\vbox{\m@th\ialign{##\crcr\noalign{\kern3\p@}
      \downbracefill\crcr\noalign{\kern3\p@\nointerlineskip}
      $\hfil\displaystyle{#1}\hfil$\crcr}}}\limits}
 \def\underbrace#1{\mathop{\vtop{\m@th\ialign{##\crcr
      $\hfil\displaystyle{#1}\hfil$\crcr\noalign{\kern3\p@\nointerlineskip}
      \upbracefill\crcr\noalign{\kern3\p@}}}}\limits}

  \let\@ldmatrix\matrix
  \let\end@ldmatrix\endmatrix
  \def\matrix{\sam\@ldmatrix}
  \def\endmatrix{\end@ldmatrix\mas}
  \let\@ldgather\gather
  \let\end@ldgather\endgather
  \def\gather{\sam\@ldgather}
  \def\endgather{\end@ldgather\mas}
  \let\@ldalign\align
  \let\end@ldalign\endalign
  \def\align{\sam\@ldalign}
  \def\endalign{\end@ldalign\mas}
  \let\@ldaligned\aligned
  \let\end@ldaligned\endaligned
  \def\aligned{\sam\@ldaligned}
  \def\endaligned{\end@ldaligned\mas}
  \let\@ldtag\tag
  \def\tag{\sam\@ldtag}
   %

   \let\MinCDArrowWidth\minCDaw@




\newskip\insertskipamount\newskip\inserthardskipamount
\insertskipamount 6pt plus2pt 
\inserthardskipamount 6pt
\def\insertskip{\vskip\insertskipamount}
\newcount\SplitTest
\def\SetSplitTest{\SplitTest\insertpenalties
  \insert\topins{\floatingpenalty1}%
  \advance\SplitTest-\insertpenalties}
\def\midinsert{\par
 \SaveLastSkip\penalty-150\SetSplitTest\RestoreLastSkip
 \ifnum\SplitTest=-1
  \@midfalse\p@gefalse\else\@midtrue\fi\@ins}
\def\@ins{\par\begingroup\setbox\z@\vbox\bgroup%
  \vglue\inserthardskipamount}
\def\endinsert{\egroup 
  \if@mid \dimen@\ht\z@ \advance\dimen@\dp\z@
    \advance\dimen@\insertskipamount
    \advance\dimen@\pagetotal\advance\dimen@-\pageshrink
    \ifdim\dimen@>\pagegoal\@midfalse\p@gefalse\fi\fi
  \if@mid%
    \ifdim\lastskip<\insertskipamount\removelastskip\insertskip\fi
    \nointerlineskip\box\z@\penalty-200\insertskip
  \else%
    \SaveLastSkip
    \insert\topins{\penalty100 
    \splittopskip\z@skip
    \splitmaxdepth\maxdimen \floatingpenalty\z@
    \ifp@ge \dimen@\dp\z@
    \vbox to\vsize{\unvbox\z@\kern-\dimen@}
    \else \box\z@\nobreak\insertskip\fi}
    \RestoreLastSkip
   \fi\endgroup}


  \newcount\notenumber
  
  \def\note{\advance\notenumber by 1
    \footnote{\the\notenumber)}}

  \newbox\footbox

  \def\footnote#1{\let\@sf\empty
    \ifhmode\edef\@sf{\spacefactor\the\spacefactor}\/\fi
    \sam${}^{\fam0 #1}$\@sf\vfootnote{#1}}%

  \def\vfootnote#1{\insert\footins\bgroup
     \interlinepenalty100 \splittopskip=1pt
     \floatingpenalty=20000
     \leftskip=0pt\rightskip=0pt%
     \parindent=.3em
     \Smallfonts\rm
     \FootItem@{#1}
     \futurelet\next\fo@t}

  \def\FootItem@#1{\par\hangafter1\hangindent=\FootHang
     \setbox0=\hbox{\ignorespaces#1\unskip}%
     \dimen0=.4em\SetOverhang@
     \noindent\rlap{\box0}\kern\Overhang\ignorespaces}


  \def\fo@t{\ifcat\bgroup\noexpand\next \let\next\f@@t
    \else\let\next\f@t\fi \next}
  \def\f@@t{\bgroup\aftergroup\@foot\let\next}
  \def\f@t#1{\baselineskip=10pt\lineskip=1pt
            \lineskiplimit=0pt #1\@foot}%
  \def\@foot{
        \hbox{\vrule height0pt depth5pt width0pt}
        \egroup}
  \skip\footins=12 pt plus 0pt minus 0pt 
  \count\footins=1000 
  \dimen\footins=8in 



 \def\osumess#1{\EdSpider{\immediate\write16{Line \the\inputlineno: #1}}}%
 \def\HideEdStuff{\gdef\EdSpider##1{}}

 \font\BigSym=cmmi10 scaled \magstep 4

 \def\change{\InLMargin{\hbox{\BigSym \char63\kern10pt}}}

 \def\beginchange{\InLMargin{\hbox{\sam\twelvepoint$\heartsuit$\kern10pt}}}

 \def\endchange{\InLMargin{\hbox{\sam\twelvepoint$\spadesuit$\kern10pt}}}

 \def\InLMargin#1{\strut\vadjust{%
     \kern-\strutdepth
     \vtop to \strutdepth{%
         \baselineskip\strutdepth
         \llap{\sam$\smash{\hbox{\EdSpider{#1}}}$}\null}}}

 \def\strutdepth{\dp\strutbox}
 \def\strutheight{\ht\strutbox}

 \def\NoteInRMargin#1{\strut\vadjust{%
     \kern-1.001\strutdepth
     \vtop to \strutdepth{%
       \baselineskip\strutdepth
       \vss\rlap{\ninepoint\unskip\hskip\hsize
         \vtop to 0pt{%
           \hsize=16em\hfuzz=\hsize
           \leftskip=10pt%
           \rightskip=0pt plus 10000pt%
           \baselineskip=9.8pt\lineskip=.2pt%
           \let\\\break
           \noindent\EdSpider{#1}\vss}%
                \kern10pt}\hbox{}}
       }}

 \def\ednote#1{\NoteInRMargin{\tentt #1}}

 \def\cbar{\InLMargin{%
      \dimen0=\strutdepth\advance\dimen0 by \lineskip
      \vrule width 3pt
      height \strutheight depth \dimen0 \kern
      3pt}}

 \def\ccbar{\InLMargin{%
      \dimen0=2\strutdepth\advance\dimen0 by 2\lineskip
      \vrule width 3pt
        height 3\strutheight depth \dimen0 \kern
      3pt}}

 \newinsert\TRMargIns
 \dimen\TRMargIns=\maxdimen

  \def\Ednote#1{\insert\TRMargIns{%
       \vbox to 0pt{\hsize=140pt\hfuzz=\hsize
           \leftskip=6pt%
           \rightskip=0pt plus 10000pt%
           \baselineskip=9.8pt\lineskip=.2pt%
           \let\\\break
           \SetPageRemainder
           \vglue540pt\vglue-\PageRemainder
           \noindent\EdSpider{\tentt #1}\vss}%
       \smallskip}}

 \def\KillEdStuff{\def\ednote##1{}\def\Ednote##1{}%
      \let\change\relax\let\beginchange\relax\let\endchange\relax
       \let\cbar\relax\let\ccbar\relax}


  \topskip=12pt
  \newskip\StdBaselineskip 
  \StdBaselineskip 12pt
  \lineskip=1.1pt
  \lineskiplimit=.8pt
  \widowpenalty=10000 
  \clubpenalty=10000  
  \abovedisplayskip=6pt plus 1pt minus 1pt
  \abovedisplayshortskip=3pt plus 1.5pt
  \belowdisplayskip=6pt plus 1pt minus 1pt
  \belowdisplayshortskip=5pt plus 1pt minus 1pt
  \hfuzz=1.5pt   

  \def\StdPretolerance{100}
  \tolerance=\StdPretolerance

  \newdimen\StdMathsurround
  \StdMathsurround=1.5pt 
  \mathsurround=\StdMathsurround
  \Mas                   

   \def\prose{\relax\hbox{\kern.6\StdMathsurround}}
  
  \def\StdParskip{0pt}    
  \parskip=\StdParskip
  \parindent=0.5cm
 

  \def\Times{ptmr  } 
  \def\TimesI{ptmri  } 
  \def\TimesB{ptmb  }
  \def\TimesBI{ptmbi  }
  \def\HelveticaN{phvrrn }

  =\Times at 10bp
  =\TimesB at 10bp
  \font\tenit=\TimesI at 10bp
  =\TimesBI at 10bp

  \font\tenmrm=cmr10  


    =\Times at 9bp 
    \font\nineit=\TimesI at 9bp 
    =\TimesB at 9bp 
    =\TimesBI at 9bp 

    =\HelveticaN at 9bp 


  =\Times at 12bp
  \font\twelveit=\TimesI at 12bp
  =\TimesB at 12bp


  \font\titleit=\TimesI at 14.4bp
  =\TimesB at 14.4bp

 \SetAuthorHead{AuthorHead} 
 \SetTitleHead{TitleHead}  


  \def\lBr{\raise.125ex\hbox{[\kern.1125ex}}
  \def\rBr{\raise.125ex\hbox{\kern.1125ex]}}

 \setbox\footbox=\hbox{\Smallfonts 2)~}



  \bgroup
  \catcode`\@=11 
  \gdef\itSpacing{%
     \xspaceskip=.31em plus.1em minus.05em \sfcode `f=2001
     \itWarning@\let\itWarning@\itWarning@@}
  \gdef\itSpacingOff{%
     \xspaceskip=0pt \sfcode `f=1000
     \let\itWarning@\relax}
   \global\let\itWarning@\relax
  \gdef\itWarning@@{\errmessage{%
  Special italic spacing already in force
  (you have probably omitted an ``endth'').
  See itSpacing macro in osuPSfnt.sty
         }}
  \egroup

 \fontdimen1\titlebf=0.0pt
 \fontdimen2\titlebf=3.6135pt
 \fontdimen3\titlebf=2.8908pt
 \fontdimen4\titlebf=1.44539pt
 \fontdimen5\titlebf=6.64882pt
 \fontdimen6\titlebf=14.45398pt
 \fontdimen7\titlebf=1.60439pt

 \fontdimen1\tenbi=0.26794pt
 \fontdimen2\tenbi=2.50937pt
 \fontdimen3\tenbi=2.00749pt
 \fontdimen4\tenbi=1.00374pt
 \fontdimen5\tenbi=4.59717pt
 \fontdimen6\tenbi=10.03749pt
 \fontdimen7\tenbi=1.11415pt

 \fontdimen1\twelverm=0.0pt
 \fontdimen2\twelverm=3.01125pt
 \fontdimen3\twelverm=2.409pt
 \fontdimen4\twelverm=1.2045pt
 \fontdimen5\twelverm=5.39615pt
 \fontdimen6\twelverm=12.045pt
 \fontdimen7\twelverm=1.33699pt

 \fontdimen1\twelveit=0.27731pt
 \fontdimen2\twelveit=3.01125pt
 \fontdimen3\twelveit=2.409pt
 \fontdimen4\twelveit=1.2045pt
 \fontdimen5\twelveit=5.37207pt
 \fontdimen6\twelveit=12.045pt
 \fontdimen7\twelveit=1.33699pt

 \fontdimen1\twelvebf=0.0pt
 \fontdimen2\twelvebf=3.01125pt
 \fontdimen3\twelvebf=2.409pt
 \fontdimen4\twelvebf=1.2045pt
 \fontdimen5\twelvebf=5.5407pt
 \fontdimen6\twelvebf=12.045pt
 \fontdimen7\twelvebf=1.33699pt

 \fontdimen1\tenrm=0.0pt
 \fontdimen2\tenrm=2.50937pt
 \fontdimen3\tenrm=2.00749pt
 \fontdimen4\tenrm=1.00374pt
 \fontdimen5\tenrm=4.49678pt
 \fontdimen6\tenrm=10.03749pt
 \fontdimen7\tenrm=1.11415pt

 \fontdimen1\tenit=0.27731pt
 \fontdimen2\tenit=2.50937pt
 \fontdimen3\tenit=2.00749pt
 \fontdimen4\tenit=1.00374pt
 \fontdimen5\tenit=4.47672pt
 \fontdimen6\tenit=10.03749pt
 \fontdimen7\tenit=1.11415pt

 \fontdimen1\tenbf=0.0pt
 \fontdimen2\tenbf=2.50937pt
 \fontdimen3\tenbf=2.00749pt
 \fontdimen4\tenbf=1.00374pt
 \fontdimen5\tenbf=4.61723pt
 \fontdimen6\tenbf=10.03749pt
 \fontdimen7\tenbf=1.11415pt

 \fontdimen1\ninerm=0.0pt
 \fontdimen2\ninerm=2.25842pt
 \fontdimen3\ninerm=1.80673pt
 \fontdimen4\ninerm=0.90337pt
 \fontdimen5\ninerm=4.0471pt
 \fontdimen6\ninerm=9.03374pt
 \fontdimen7\ninerm=1.00273pt

 \fontdimen1\nineit=0.27731pt
 \fontdimen2\nineit=2.25842pt
 \fontdimen3\nineit=1.80673pt
 \fontdimen4\nineit=0.90337pt
 \fontdimen5\nineit=4.02904pt
 \fontdimen6\nineit=9.03374pt
 \fontdimen7\nineit=1.00273pt

 \fontdimen1\ninebf=0.0pt
 \fontdimen2\ninebf=2.25842pt
 \fontdimen3\ninebf=1.80673pt
 \fontdimen4\ninebf=0.90337pt
 \fontdimen5\ninebf=4.15552pt
 \fontdimen6\ninebf=9.03374pt
 \fontdimen7\ninebf=1.00273pt


 \newcount\MaxSpaceFactor
 \MaxSpaceFactor=3000 

 \def\ItemStyle{\rm}
 \def\NrStyle{\rm}
 \def\ItemItemStyle{\rm}

 \MaxItemTag{(iii)}
 \MaxItemItemTag{(iii)}
 \MaxNrTag{(2)}
 \MaxFootTag{2)}
 \def\ReferenceHang{30pt}

 \catcode`\@=\active


\loadbold

=\Times  
=\Times scaled750
=\Times scaled650
\font\rms=\Times scaled 920 

=\TimesBI scaled 860
=\TimesI scaled 860

\textfont0=\rrm  
\scriptfont0=\erm 
\scriptscriptfont0=\srm

\def\Augment#1#2{%
    \toks0\expandafter{#1}\toks2{#2}%
    \edef#1{\the\toks0\the\toks2}}

 \font\twelverma=\Times  scaled 1200
 \font\tenrma=\Times  scaled 1000
 \font\ninerma=\Times scaled 920
 =\Times scaled 840
 \font\sevenrma=\Times scaled 760
 =\Times scaled 680
 \font\fiverma=\Times scaled 600

 \Augment\tenpoint{%
  \textfont0=\tenrma  \scriptfont0=\sevenrma  
  \scriptscriptfont0=\fiverma  }

 \Augment\ninepoint{%
  \textfont0=\ninerma  \scriptfont0=\sevenrma 
  \scriptscriptfont0=\fiverma}

 \Augment\twelvepoint{%
  \textfont0=\twelverma  \scriptfont0=\ninerma  
  \scriptscriptfont0=\sevenrma}

\mathsurround=1pt
\hsize=13.45truecm
\vsize=19.5truecm
\hoffset=1.25truecm
\voffset=2truecm
\advance\baselineskip by 2pt

\predefine\til{\~}
\def\~#1{\relax\ifmmode\widetilde{#1}\else\til{#1}\fi}

\redefine \le{\leqslant}
\redefine \ge{\geqslant}
\define \wt#1{\mathaccent"0365{#1}}
\define \wh#1{\mathaccent"0362{#1}}

\define \iss{\,\Mathaccent{\raise -.8 ex\hbox{$\widetilde{}$\kern.1em}}\rightarrow\,}

\define \tpp{\mathop{\fam0 top}}
\define \ab{\mathop{\fam0 ab}}

\define \insep{\mathop{\fam0 ins}}
\define \alg{\mathop{\fam0 alg}}

\define \chr{\mathop{\fam0 char}\,}

\define \Gal{\mathop{\fam0 Gal}}

\define \Fil{\operatorname{\fam0 fil}}

\Mas
\HideEdStuff
\rm 
 

\def\issn{{\nineit ISSN 1464-8997 (on line) 1464-8989 (printed)}}

\def\gtp{{\nineit Published 10 December 2000: \ \copyright\ Geometry \& 
Topology Publications}}

\def\gtv3{{\nineit Geometry \& Topology Monographs, Volume 3 (2000) --
Invitation to higher local fields}}


\def\lione
{{\rms Geometry \& Topology Monographs}}

\def \litwo{{\rms Volume 3: Invitation to higher local fields
}} 

\def\tinfo #1.#2.#3-#4
{{
\noindent  {\lione} \hfill 
\par 
\vskip-1.5pt
\noindent {\litwo} \hfill
\par 
\vskip-1,5pt
\noindent {\rms Part #1, section #2, pages #3--#4} \hfill
\vskip24pt 
}}

\def\tinfos #1.#2.#3-#4
{{
\noindent  {\lione} \hfill 
\par 
\vskip-1.5pt
\noindent {\litwo} \hfill
\par 
\vskip-1.5pt
\noindent {\rms Pages #3--#4} \hfill
\vskip24pt 
}}

\def\tinfoi #1
{{
\noindent  {\lione} \hfill 
\par 
\vskip-1.5pt
\noindent {\litwo} \hfill
\par 
\vskip-1.5pt
\noindent {\rms Pages iii--xi: Introduction and contents} \hfill
\vskip26pt 
}}


  \def\titlepagehead{\hfil}

  \newif\iftitlepage\titlepagefalse
  \newif\ifblankpage\blankpagefalse
  \def\makeheadline{
     \ifblankpage{}\else%
     \iftitlepage
\vbox{\line{\vbox to 8.5pt{}
\ninerm
\copy\HLinebox \hfill
\hglue5mm\ninebf\folio 
\titlepagehead}}%
      \else
\vbox{\ifodd\pageno\rightheadline\else\leftheadline\fi}%
      \fi\vskip 12pt\fi}%
     \def\rightheadline{\line{\vbox to 8.5pt{}%
      \ninerm
\copy\TitleBox \hfill
\hglue5mm\ninebf\folio}}%
     \def\leftheadline{\line{\vbox to 8.5pt{}%
        \unskip\ninerm\unskip\ninebf\folio\hglue5mm
 \hfill \copy\AuthorBox
}}

 \footline={\ifblankpage{}\else
\iftitlepage\ninepoint\sam\hfill
\line{\vbox to 8.5pt{}
\copy\TFLinebox
\hfill
\hglue5mm 
}
            \else
\ninepoint\sam\hfill
\line{\vbox to 8.5pt{}
\copy\FLinebox
\hfill 
\hglue5mm
}
\hfil\fi\global\titlepagefalse\fi}

\def\blankpage{{\blankpagetrue\noindent\hbox to 10pt{\hss}\vfill
\pagebreak}}

\tenpoint\rm 
 

\pageno=143

\tinfo I.17.143-150

\SetTFLinebox{\gtp }
\SetFLinebox{\gtv3 }
\SetHLinebox{\issn}

\H 17. An approach to higher ramification theory

Igor Zhukov

\SetAuthorHead{I. Zhukov}
\SetTitleHead{Part I. Section 17. An approach to higher ramification theory
 \qquad\qquad}

We use the notation of sections~1 and 10.

\HH 17.0. Approach of Hyodo and Fesenko

Let $K$ be an $n$-dimensional local field,
$L/K$ a finite abelian extension.
Define a filtration on $\Gal(L/K)$ (cf. \cite{H}, \cite{F, sect. 4}) by
$$
\Gal(L/K)^{\bold i}=\Upsilon_{L/K}^{-1} (U_{{\bold i} }K_n^{\tpp}(K)+N_{L/K}K_n^{\tpp}(L)/
N_{L/K}K_n^{\tpp}(L)),\quad {\bold i}\in \Bbb Z^n_+,
$$
where $U_{\bold i} K_n^{\tpp}(K)=\{U_{\bold i}\}\cdot K_{n-1}^{\tpp}(K)$,
$U_{\bold i}=1+P_K({\bold i})$, 
$$
\Upsilon_{L/K}^{-1}\colon K_n^{\tpp}(K)/N_{L/K}K_n^{\tpp}(L)\iss \Gal(L/K)
$$
is the reciprocity map.

Then for a subextension $M/K$ of $L/K$
$$\Gal(M/K)^{{\bold i}} =\Gal(L/K)^{{\bold i}} \Gal(L/M)/\Gal(L/M)$$
which is a higher dimensional analogue of Herbrand's theorem.
However, if one defines a generalization of 
the Hasse--Herbrand function and  lower ramification filtration,
then for $n>1$ the lower filtration on a subgroup does not coincide with the induced filtration in general. 

Below we shall give another construction of the ramification filtration
of $L/K$ in the two-dimen\-si\-o\-nal case;  details can be found
in \cite{Z}, see also \cite{KZ}. This construction can be considered as a development
of an approach by K.~Kato and T.~Saito in \cite{KS}.

\df Definition

Let $K$ be a complete discrete valuation field
with residue field $k_K$ of characteristic $p$.
A finite extension $L/K$ is called {\it ferociously ramified}
if $|L:K|=|k_L:k_K|_{\insep}$.
\enddf

In addition to the nice ramification theory for 
totally ramified extensions,
there is a nice ramification theory for ferociously 
ramified extensions $L/K$ such that $k_L/k_K$ is generated by one element; 
the reason is that in both cases the ring extension $\Cal O_L/\Cal O_K$ is monogenic,
i.e., generated by one element, see section 18. 

\HH 17.1. Almost constant extensions 

Everywhere below  $K$ is a complete discrete valuation field
with residue field $k_K$ of characteristic $p$
such that  $|k_K:k_K^p|=p$.
For instance, $K$ can be a two-dimensional local field, or $K=\Bbb F_q(X_1)((X_2))$
or the quotient field of the completion of $\Bbb Z_p[T]_{(p)}$
with respect to the $p$-adic topology.

\df Definition

For the field $K$ define a base (sub)field $B$ as

$B=\Bbb Q_p\subset K$ if $\chr (K)=0$,

$B=\Bbb F_p((\rho))\subset K$
if $\chr(K)=p$, where $\rho$ is an element of $K$ with $v_K(\rho)>0$.

Denote by $k_0$ the completion of $B(\Cal R_K)$ inside $K$. 
Put $k=k_0^{\alg}\cap K$.

The subfield $k$ is a maximal complete subfield of $K$ 
with perfect residue field. 
It is called a {\it constant subfield} of $K$.
A constant subfield is defined canonically if $\chr (K)=0$.
Until the end
 of section 17 we assume that $B$ (and, therefore, $k$) is fixed.

By $v$ we denote the valuation ${K^{\alg}}^*\to\Bbb Q$ normalized so 
that $v(B^*)=\Bbb Z$.

\enddf

\eg Example

If $K=F\{\!\{T\}\!\}$ where $F$ is a mixed characteristic complete discrete valuation field with perfect residue field, then $k=F$. 
\endeg

\df Definition

An extension $L/K$ is said to be {\it constant}
if there is an algebraic extension $l/k$ such that
$L=Kl$.

An extension $L/K$ is said to be {\it almost constant}
if $L\subset L_1L_2$ for a constant extension
$L_1/K$ and an unramified extension $L_2/K$.

A field $K$ is said to be {\it standard}, if $e(K|k)=1$, and {\it almost standard},
if some finite unramified extension of $K$ is a standard field.
\enddf

\th Epp's theorem on
elimination of wild ramification

{\tenrm (\cite{E}, \cite{KZ})}
Let $L$ be a finite extension of $K$.
Then there is a finite extension $k'$ of a constant subfield $k$ of $K$
such that $e(Lk'|Kk')=1$.
\endth

\th Corollary 

There exists a finite constant extension of $K$ which is a standard field.
\endth

\pf Proof

See the proof of the Classification Theorem in { 1.1}.
\endpf

\th Lemma

The class of constant (almost constant) extensions
is closed with respect to
taking compositums and subextensions. If $L/K$ and $M/L$ are almost constant
then $M/K$ is almost constant as well.
\endth

\df Definition

Denote by $L_c$ the maximal almost constant subextension
of $K$ in $L$.
\enddf

\rk Properties

\Roster

\Item{(1)} 
Every tamely ramified extension is almost constant. In other words,
the (first) ramification subfield in $L/K$ is a subfield
of $L_c$.

\Item{(2)} If $L/K$ is normal then
$L_c/K$ is normal.

\Item{(3)} There is an unramified extension
$L_0'$ of $L_0$ such that
$L_cL_0'/L_0$ is a constant extension.

\Item{(4)} There is a constant extension $L_c'/L_c$ such that
$LL_c'/L_c'$ is ferociously ramified and $L_c'\cap L=L_c$.
This follows  immediately from Epp's theorem.
\endRoster 
\endrk

The principal idea of the proposed approach to ramification theory
is to split $L/K$ into a tower of three extensions:
$L_0/K$, $L_c/L_0$, $L/L_c$, where $L_0$ is the inertia subfield
in $L/K$. The ramification filtration for $\Gal(L_c/L_0)$ reflects that
for the corresponding extensions of constants subfields. Next,
to construct the ramification filtration for $\Gal(L/L_c)$, one reduces
to the case of ferociously ramified extensions by means of Epp's theorem.
(In the case of higher local fields one can also construct
a filtration on $\Gal(L_0/K)$ by lifting that for the first residue fields.)

Now we give precise definitions.

\HH 17.2. Lower and upper ramification filtrations

Keep the assumption of the previous subsection.
Put
$$\Cal A=\{
-1,0\}\cup \{({\goth c},s):0<s\in\Bbb Z\}\cup \{({\goth i},r): 0<r\in\Bbb Q\}.
$$
This set is
linearly ordered as follows:
$$ \gather -1<0<({\goth c},i)<({\goth i},j) \text{ for any }i,j;
 \\  ({\goth c},i)<({\goth c},j) \text{ for any }i<j;
\\ ({\goth i},i)<({\goth i},j) \text{ for any }i<j.
\endgather
$$

\df Definition

Let  $G=\Gal(L/K)$. For any
$\alpha\in\Cal A$ we  define a subgroup $G_\alpha$ in $G$.

Put $G_{-1}=G$, and denote by $G_0$ the inertia subgroup
in $G$, i.e.,
$$
G_0=\{g\in G : v(g(a)-a)>0\text{ for all }a\in\Cal O_L\}.
$$

Let $L_c/K$ be constant, and let it contain no unramified subextensions.
Then define 
$$G_{{\goth c},i}={\text{\tenrm pr}}^{-1}(\Gal(l/k)_i)$$ 
where
$l$ and $k$ are the constant subfields in $L$ and $K$ respectively,
$${\text{\tenrm pr}}\:\Gal(L/K)\to\Gal(l/k)=\Gal(l/k)_0$$ is the natural projection
and $\Gal(l/k)_i$ are the classical ramification subgroups.
In the general case take an unramified extension $K'/K$
such that $K'L/K'$ is constant and contains no unramified subextensions, 
and put
$G_{{\goth c},i}=\Gal(K'L/K')_{{\goth c},i}$.

Finally, define $G_{{\goth i},i}$, $i>0$.
Assume that $L_c$ is standard and $L/L_c$ is ferociously ramified.
Let $t\in \Cal O_L$, $\Overline{t}\notin k_L^p$.
Define
$$
G_{{\goth i},i}=\{g\in G: v(g(t)-t)\ge i\}
$$
for all $i>0$. 

In the general case choose a finite extension $l'/l$ such that $l'L_c$ is standard and
$e(l'L|l'L_{c})=1$. Then it is clear that
$\Gal(l'L/l'L_c)=\Gal(L/L_c)$,
and $l'L/l'L_c$ is ferociously ramified. Define
$$
G_{{\goth i},i}=\Gal(l'L/l'L_c)_{{\goth i},i}
$$
for all $i>0$.
\enddf

\th Proposition

For a finite Galois extension $L/K$ the lower filtration
$\{\Gal(L/K)_\alpha\}_{\alpha\in \Cal A}$ is well defined.
\endth

\df Definition

Define a generalization  $h_{L/K}\colon \Cal A\to\Cal A$
of the Hasse--Herbrand function. 
First, we define  $$\Phi_{L/K}\:\Cal A\to\Cal A$$ as follows:
$$
\aligned
\Phi_{L/K}(\alpha)&=\alpha  \quad\text{ for } \alpha=-1,0;
 \\
\Phi_{L/K}(({\goth c},i))&=
 \biggl({\goth c},\frac1{e(L|K)}\int_0^i|\Gal(L_c/K)_{{\goth c},t}|dt\biggr)
\quad \text{ for all } i>0;
 \\
\Phi_{L/K}(({\goth i},i))&=
\biggl({\goth i},\int_0^i|\Gal(L/K)_{{\goth i},t}|dt\biggr)
\quad \text{ for all } i>0.
\endaligned 
$$
It is easy to see that $\Phi_{L/K}$ is bijective and increasing,
and we introduce $$h_{L/K}=\Psi_{L/K}=\Phi_{L/K}^{-1}.$$

Define the upper filtration
$\Gal(L/K)^\alpha=\Gal(L/K)_{h_{L/K}(\alpha)}$.
\enddf 

All standard formulas for intermediate extensions take place;
in particular, for a normal subgroup $H$ in $G$ we have
$H_\alpha=H\cap G_\alpha$ and $(G/H)^\alpha=G^\alpha H/H$. The latter relation
enables one to introduce the upper filtration for an infinite
Galois extension as well.

\rk Remark

The filtrations do depend on the choice of a constant subfield
(in characte\-ris\-tic~$p$).
\endrk

\eg Example

Let $K=\Bbb F_p((t))((\pi))$.
Choose $k=B=\Bbb F_p((\pi))$ as a constant subfield.
Let $L=K(b)$, $b^p-b=a\in K$.
Then
\Roster

\Item{} if $a=\pi^{-i}$, $i$ prime to $p$, then the ramification break
of $\Gal(L/K)$ is $(\goth c,i)$;

\Item{} if $a=\pi^{-pi}t$, $i$ prime to $p$, then the ramification break
of $\Gal(L/K)$ is $(\goth i,i)$;

\Item{} if $a=\pi^{-i}t$, $i$ prime to $p$, then the ramification break
of $\Gal(L/K)$ is $(\goth i,i/p)$;

\Item{} if $a=\pi^{-i}t^p$, $i$ prime to $p$, then the ramification break
of $\Gal(L/K)$ is $(\goth i,i/p^2)$.
\endRoster
\endeg

\rk Remark

A dual filtration on $K/\wp(K)$ is computed in the final version of \cite{Z}, see 
also \cite{KZ}.
\endrk

\HH 17.3. Refinement for a two-dimensional local field

Let $K$ be a two-dimensional local field with $\chr (k_K)=p$, and
let $k$ be the constant subfield of $K$. Denote by
$$
{\bold v}=(v_1,v_2)\colon (K^{\alg})^*\to \Bbb Q\times \Bbb Q
$$
the extension of the rank 2 valuation of $K$, which is normalized so that:

$\bullet$ $v_2(a)=v(a)$ for all $a\in K^*$,

$\bullet$ $v_1(u)=w({\Overline u})$ for all $u\in U_{K^{\alg}}$, where $w$ is a non-normalized extension
of $v_{k_K}$ on $k_K^{\alg}$, and $\Overline u$ is the residue of $u$,

 $\bullet$ ${\bold v}(c)=(0,e(k|B)^{-1}v_k(c))$ for all $c\in k$.

It can be easily shown that ${\bold v}$ is uniquely determined by these conditions,
and the value group of ${\bold v}|_{K^*}$ is isomorphic to $\Bbb Z\times \Bbb Z$.

Next, we introduce the index set
$$
\Cal A_2=\Cal A\cup\Bbb Q_+^2=\Cal A\cup\{(i_1,i_2) : i_1,i_2\in\Bbb Q,i_2>0\}
$$
and extend the ordering of $\Cal A$ onto $\Cal A_2$ assuming
$$
(\goth i,i_2)<(i_1,i_2)<(i_1',i_2)<(\goth i,i_2')
$$
for all $i_2<i_2'$, $i_1<i_1'$.

Now we can define $G_{i_1,i_2}$, where $G$ is the Galois
group of a given finite Galois extension $L/K$. Assume first that $L_c$ is standard and
$L/L_c$ is ferociously ramified.
Let $t\in \Cal O_L$, $\bar t\notin k_L^p$ (e.g., a first local parameter of $L$).
We define
$$
G_{i_1,i_2}=\bigl\{g\in G : 
{\bold v}\bigl(t^{-1}g(t)-1\bigr)\ge(i_1,i_2)\bigr\}
$$
for $i_1,i_2\in\Bbb Q$, $i_2>0$. In the general case we choose $l'/l$ 
($l$ is the constant subfield of both $L$ and $L_c$)
such that $l'L_c$ is standard and $l'L/l'L_c$ is ferociously ramified and put
$$
G_{i_1,i_2}=\Gal(l'L/l'L_c)_{i_1,i_2}.
$$
We obtain a well defined lower filtration $(G_\alpha)_{\alpha\in\Cal A_2}$ on $G=\Gal(L/K)$.

In a similar way to 17.2, one constructs the Hasse--Herbrand functions\hfill
\break
$\Phi_{2,L/K}\:\Cal A_2\to\Cal A_2$ and $\Psi_{2,L/K}=\Phi_{2,L/K}^{-1}$ which extend
$\Phi$ and $\Psi$ respectively. Namely,
$$
\Phi_{2,L/K}((i_1,i_2))=\int_{(0,0)}^{(i_1,i_2)}|\Gal(L/K)_{t}|dt.
$$

These functions have usual properties of the Hasse--Herbrand functions $\varphi$ and $h=\psi$,
and one can introduce an $\Cal A_2$-indexed upper filtration on any finite
or infinite Galois group $G$.

\HH 17.4. Filtration on $K^{\tpp}(K)$

In the case of a two-dimensional local field $K$ the upper ramification filtration
for $K^{\ab}/K$ determines a compatible filtration on $K_2^{\tpp}(K)$. In the
case where $\chr (K)=p$ this filtration has an explicit description given below.

From now on, let $K$ be a two-dimensional local field of prime characteristic $p$ over a quasi-finite field,
and $k$ the constant subfield of $K$.
Introduce ${\bold v}$ as in {17.3}. Let $\pi_k$ be a prime of $k$.

For all $\alpha\in\Bbb Q_+^2$ introduce subgroups
$$
\aligned
&Q_\alpha =\{\,\{\pi_k,u\}\, : \,u\in K, {\bold v}(u-1)\ge\alpha\,\}\subset VK_2^{\tpp} (K); \\
&Q_\alpha^{(n)}=\{a\in K_2^{\tpp} (K) \, : \, p^na\in Q_\alpha\}; \\
&S_\alpha ={\text{\rm Cl}\,} \bigcup_{n\ge0} Q_{p^n\alpha}^{(n)}.
\endaligned
$$
For a subgroup $A$,  ${\text{\rm Cl}\,} A$ denotes the intersection of all open subgroups containing $A$.

The subgroups $S_\alpha$ constitute the heart of the ramification filtration on $K_2^{\tpp}(K)$.
Their most important property is that they have nice behaviour in unramified, constant
and ferociously ramified extensions.

\th Proposition 1

Suppose that $K$ satisfies the following property.

\Roster
\Item{\rm (*)}  The extension of constant subfields
in any finite unramified extension of $K$ is also unramified.  
\endRoster

Let $L/K$ be either an unramified or a constant totally ramified extension, $\alpha\in \Bbb Q^2_+$. Then we have
$N_{L/K}S_{\alpha,L}=S_{\alpha,K}$.
\endth

\th Proposition 2

Let $K$ be standard, $L/K$ a cyclic ferociously ramified extension of degree $p$
with the ramification jump $h$ in lower numbering, $\alpha\in\Bbb Q^2_+$.
Then:

{\rm (1)} $N_{L/K}S_{\alpha,L}=S_{\alpha+(p-1)h,K}$, if $\alpha>h$;

{\rm (2)} $N_{L/K}S_{\alpha,L}$ is a subgroup in $S_{p\alpha,K}$ of index $p$, if $\alpha\le h$.
\endth

Now we have ingredients to define a decreasing filtration
$\{\Fil_\alpha K_2^{\tpp} (K)\}_{\alpha\in\Cal A_2}$ on $K_2^{\tpp} (K)$.
Assume first that $\~K$ satisfies the condition (*).
It follows from \cite{KZ, Th. 3.4.3} that for some purely inseparable constant extension $K'/K$
the field $K'$ is almost standard. Since $K'$ satisfies (*) and is almost standard,
it is in fact standard.

Denote
$$
\alignat 2
&\Fil_{\alpha_1,\alpha_2}K_2^{\tpp} (K)&& =S_{\alpha_1,\alpha_2}; \\
&\Fil_{{\goth i},\alpha_2}K_2^{\tpp} (K)&& ={\text{\rm Cl}\,}\bigcup\limits_{\alpha_1\in\Bbb Q}\Fil_{\alpha_1,\alpha_2} K_2^{\tpp} (K)
 \text{ for } \alpha_2\in\Bbb Q_+; \\
&T_K&&={\text{\rm Cl}\,}\bigcup\limits_{\alpha\in\Bbb Q^2_+}\Fil_{\alpha}K_2^{\tpp} (K);\\
&\Fil_{{\goth c},i}K_2^{\tpp}(K)&&=T_K+\{\,\{t,u\}\,: \,u\in k,\, v_k(u-1)\ge i\}\text{ for all $i\in\Bbb Q_+$,}\\
&&&\quad \text{ if $K=k\{\!\{t\}\!\}$ is standard};\\
&\Fil_{{\goth c},i}K_2^{\tpp}(K)&&=N_{K'/K}\Fil_{{\goth c},i}K_2^{\tpp}(K')\text{, where $K'/K$ is as above};\\
&\Fil_0K_2^{\tpp}(K)&&=U(1)K_2^{\tpp}(K)+\{t,\Cal R_K\}\text{, where } 
U(1)K_2^{\tpp}(K)=\{1+P_K(1),K^*\},\\
&&&\quad \text{ $t$ is the first local parameter};\\
&\Fil_{-1}K_2^{\tpp}(K)&&=K_2^{\tpp}(K).
\endalignat
$$

It is easy to see that
for some unramified extension  $\~K /K$ the field $\~K$ satisfies the condition (*),
and we define $\Fil_\alpha K_2^{\tpp}(K)$ as 
$N_{\~K/K}\Fil_\alpha K_2^{\tpp}(\~K)$ for all $\alpha\ge0$,
and $\Fil_{-1} K_2^{\tpp}(K)$ as $K_2^{\tpp}(K)$.
It can be shown that the filtration $\{\Fil_\alpha K_2^{\tpp}(K)\}_{\alpha\in\Cal A_2}$ is well defined.

\th Theorem 1

Let $L/K$ be a finite abelian extension, $\alpha\in\Cal A_2$. Then $N_{L/K}\Fil_\alpha K_2^{\tpp} (L)$
is a subgroup in $\Fil_{\Phi_{2,L/K}(\alpha)}K_2^{\tpp}(K)$ of index $|\Gal(L/K)_{\alpha}|$.
Furthermore,
$$
\Fil_{\Phi_{L/K}(\alpha)}K_2^{\tpp}(K)\cap N_{L/K}K_2^{\tpp} (L)=N_{L/K}\Fil_\alpha K_2^{\tpp} (L).
$$
\endth

\th Theorem 2

Let $L/K$ be a finite abelian extension, and let
$$\Upsilon_{L/K}^{-1}\:K_2^{\tpp}(K)/N_{L/K}K_2^{\tpp} (L)\to \Gal(L/K)$$ 
be the reciprocity map. Then
$$
\Upsilon_{L/K}^{-1}(\Fil_\alpha K_2^{\tpp}(K) \mod N_{L/K}K_2^{\tpp} (L))= \Gal(L/K)^\alpha
$$
for any  $\alpha\in\Cal A_2$.
\endth

\rk Remarks

1. The ramification filtration, constructed in {17.2},  does not give  information 
about the classical ramification invariants in general. Therefore, this construction can
be considered only as a provisional one.

2. The filtration on $K_2^{\tpp}(K)$ constructed in {17.4} behaves 
with
respect to the norm map much better  than the usual filtration $\{U_{{\bold i}}K_2^{\tpp}(K)\}_{{\bold i}\in\Bbb Z^n_+}$. We hope that this
filtration can be useful in the study of the structure of $K^{\tpp}$-groups.

3. In the mixed characteristic case the description of ``ramification'' filtration  on $K_2^{\tpp}(K)$ is not very nice.
However, it would be interesting to try to modify the ramification filtration
on $\Gal(L/K)$ in order to get  the filtration  on $K_2^{\tpp}(K)$ similar to that described in  {17.4}.

4. It would be interesting to compute ramification of the extensions
constructed in sections 13 and~14.
\endrk

\Bib References 

\rf{E}
H.~Epp,
Eliminating wild ramification,
Invent. Math. 19 (1973),
pp. 235--249

\rf{F}
I. Fesenko, 
Abelian local $p$-class field theory,  
 Math. Ann. 301 (1995), 561--586.

\rf{H} O. Hyodo,  {Wild ramification in the imperfect residue 
field case, } Advanced Stud. in Pure Math. {12} (1987) Galois 
Representation and Arithmetic  Algebraic Geometry, 287--314.

\rf{KS}
K. Kato and T. Saito,
Vanishing cycles, ramification of valuations and class field theory,
Duke Math. J.,
55 (1997), 629--659

\rf{KZ}
 M. V. Koroteev and I. B. Zhukov,
 Elimination of wild ramification,
 Algebra i Analiz
11 (1999), no. 6.


\rf{Z}
I. B. Zhukov,
On ramification theory in the imperfect residue
field case,
preprint of Nottingham University 98-02,
Nottingham, 1998, www.dpmms.cam.ac.uk/Algebraic-Number-Theory/0113,
to appear in Proceedings of the Luminy conference on Ramification theory
for arithmetic schemes.

\endBib

\Coordinates

Department of Mathematics and Mechanics \ 
St. Petersburg University

Bibliotechnaya pl. 2, \ 
Staryj Petergof

198904 St. Petersburg \ Russia

E-mail: igor\@zhukov.pdmi.ras.ru

\endCoordinates

\vfill
\pagebreak

\bye